\documentclass[11pt]{amsart}

\usepackage[style=numeric]{biblatex} 
\usepackage[T1]{fontenc}
\usepackage{lmodern,microtype}
\usepackage[margin=1in]{geometry}
\usepackage{amsmath,amssymb,amsthm,mathtools}
\usepackage[colorlinks=true,linkcolor=blue!50!black,citecolor=blue!50!black,urlcolor=blue!50!black]{hyperref}
\usepackage{xcolor}

\newtheorem{theorem}{Theorem}[section]
\newtheorem{proposition}[theorem]{Proposition}
\newtheorem{lemma}[theorem]{Lemma}
\newtheorem{corollary}[theorem]{Corollary}
\theoremstyle{definition}

\newtheorem{remark}[theorem]{Remark}

\newcommand{\Hess}{\operatorname{Hess}}
\newcommand{\tr}{\operatorname{tr}}
\newcommand{\Aff}{\operatorname{Aff}}
\newcommand{\Sc}{\mathcal S}

\newcommand{\Ree}{\operatorname{Re}}

\begin{document}

\title{Splitting and Mixed Deformations of\\Toric Generalized K\"ahler Potentials}
\author{Hongnian Huang \and Li Sheng}
\address{
    Hongnian Huang,
    Department of Mathematics and Statistics,
    University of New Mexico,
    Albuquerque, NM 87131, USA
}
\email{hnhuang@gmail.com}

\address{
Li Sheng,
Department of Mathematics,
Sichuan University,
Chengdu 610065, Sichuan, China
}
\email{lsheng@scu.edu.cn}

\begin{abstract}
We study splitting phenomena for extremal toric generalized K\"ahler structures on products of compact toric manifolds. We prove that every extremal symplectic potential splits when the skew-symmetric parameter \(B\) is block diagonal, without any restriction on \(A\). We also construct two types of non-splitting examples on \(\mathbb{CP}^1\times\mathbb{CP}^1\): an explicit constant-generalized-scalar-curvature structure with a splitting potential but a non-product generalized K\"ahler structure, and a deformation family of constant-generalized-scalar-curvature structures with genuinely non-splitting potentials.
\end{abstract}

\maketitle

\section{Introduction}
The product question studied here concerns canonical geometric structures. In the K\"ahler geometry, canonical K\"ahler metrics exhibit remarkable rigidity with respect to product decompositions. Apostolov and the first named author \cite{ApostolovHuang} proved a splitting result in extremal Kähler metrics on products of polarized projective manifolds, under certain assumptions involving integral Futaki invariants or automorphism groups. In the toric setting, the first named author \cite{HuangSplitting} removed these additional assumptions and proved that every extremal Kähler metric on a product of compact toric manifolds splits as a product of extremal Kähler metrics on the factors, using the short-time existence of the Calabi flow. These results naturally motivate the investigation of analogous splitting phenomena in generalized Kähler geometry. In the present paper, we establish a splitting theorem for the symplectic potentials of extremal toric generalized Kähler structures when the skew-symmetric parameter \(B\) respects the product decomposition. We also construct two types of examples on \(\mathbb{CP}^1\times\mathbb{CP}^1\) illustrating distinct failures of splitting: the first consists of extremal generalized Kähler structures whose symplectic potentials split but whose full generalized Kähler structures do not, while the second consists of extremal generalized Kähler structures whose symplectic potentials themselves fail to split.

\subsection{Generalized K\"ahler geometry}
Generalized Kähler geometry originated in theoretical physics through the work of Gates, Hull, and Roček \cite{GatesHullRocek}, who discovered the bi-Hermitian geometry underlying two-dimensional supersymmetric nonlinear sigma models. Their work showed that the target spaces of sigma models with \(N=(2,2)\) supersymmetry admit a more general geometric structure than ordinary Kähler geometry, involving two integrable complex structures compatible with a common Riemannian metric and a torsion three-form. A new geometric perspective emerged with Hitchin's introduction of generalized complex geometry \cite{Hitchin2003}, which unifies complex and symplectic geometry through structures on the generalized tangent bundle \(TM\oplus T^*M\). Building on Hitchin's framework, Gualtieri \cite{GualtieriThesis,Gualtieri} introduced generalized Kähler structures and established their equivalence with the bi-Hermitian geometry arising in the earlier physics literature.

Let \(M\) be a smooth manifold of real dimension \(2n\). The generalized tangent bundle is the direct sum
\[
\mathbb{T}M:=TM\oplus T^*M,
\]
equipped with the natural symmetric bilinear form of signature \((2n,2n)\), defined by
\[
\langle X+\xi,Y+\eta\rangle
=\frac12\bigl(\xi(Y)+\eta(X)\bigr),
\]
for \(X,Y\in TM\) and \(\xi,\eta\in T^*M\). The space of smooth sections of \(\mathbb{T}M\) is endowed with the Courant bracket
\[
[X+\xi,Y+\eta]_C
=[X,Y]+\mathcal L_X\eta-\mathcal L_Y\xi
-\frac12d\bigl(\iota_X\eta-\iota_Y\xi\bigr).
\]

A generalized almost complex structure is an endomorphism \(\mathcal J:\mathbb{T}M\to\mathbb{T}M\) satisfying
\[
\mathcal J^2=-\mathrm{Id},
\qquad
\langle\mathcal Ju,\mathcal Jv\rangle
=\langle u,v\rangle
\]
for all \(u,v\in\mathbb{T}M\). Equivalently, \(\mathcal J\) is an orthogonal complex structure on the generalized tangent bundle. It is said to be integrable if its \(+\mathrm{i}\)-eigenbundle
\[
L_{\mathcal J}
=\ker(\mathcal J-\mathrm{i}\,\mathrm{Id})
\subset\mathbb{T}M\otimes\mathbb C
\]
is involutive under the complexified Courant bracket, that is,
\[
[\Gamma(L_{\mathcal J}),\Gamma(L_{\mathcal J})]_C
\subseteq\Gamma(L_{\mathcal J}).
\]
An integrable generalized almost complex structure is called a generalized complex structure.

A generalized Kähler structure is a pair \((\mathcal J_1,\mathcal J_2)\) of commuting integrable generalized complex structures,
\[
[\mathcal J_1,\mathcal J_2]=0,
\]
such that
\[
\mathcal G:=-\mathcal J_1\mathcal J_2
\]
defines a positive generalized metric. More precisely, \(\mathcal G\) is an orthogonal, self-adjoint involution satisfying
\[
\mathcal G^2=\mathrm{Id},
\qquad
\langle\mathcal Gu,v\rangle
=\langle u,\mathcal Gv\rangle,
\]
and the bilinear form
\[
(u,v)_{\mathcal G}:=\langle\mathcal Gu,v\rangle
\]
is positive definite. Consequently, the generalized tangent bundle admits an orthogonal decomposition
\[
\mathbb{T}M=V_+\oplus V_-,
\]
where \(V_\pm\) are the \(\pm1\)-eigenbundles of \(\mathcal G\). The natural pairing is positive definite on \(V_+\) and negative definite on \(V_-\).

A generalized metric is equivalently determined by a Riemannian metric \(g\) and a real two-form \(\beta\). In the chosen splitting of the Courant algebroid, its eigenbundles are given by
\[
V_\pm
=\left\{
X+(\beta\pm g)(X):X\in TM
\right\},
\]
where \(g\) and \(\beta\) are regarded as bundle maps from \(TM\) to \(T^*M\). Thus, the positivity condition on \(\mathcal G\) recovers an ordinary Riemannian metric on \(M\), while the commuting generalized complex structures encode additional complex-geometric data. 

Gualtieri's bi-Hermitian correspondence \cite{GualtieriThesis,Gualtieri} establishes the equivalence between generalized Kähler structures and the bi-Hermitian geometry originally arising in supersymmetric sigma models. More precisely, upon choosing an isotropic splitting of the underlying exact Courant algebroid, a generalized K\"ahler structure can be described by a quadruple \((g,\beta,I,J)\), where \(g\) is a Riemannian metric, \(\beta\) is a real two-form, and (I,J) are two integrable complex structures compatible with \(g\), namely,
\[
g(IX,IY)=g(JX,JY)=g(X,Y)
\]
for all \(X,Y\in TM\). The corresponding fundamental two-forms are defined by
\[
\omega_I(X,Y)=g(IX,Y),
\qquad
\omega_J(X,Y)=g(JX,Y).
\]
In the untwisted setting, the integrability conditions for the generalized K\"ahler structure are equivalent to the torsion identities
\[
\boxed{
-d_I^c\omega_I=d\beta=d_J^c\omega_J.
}
\]
In particular, the two Hermitian structures \((g,I)\) and \((g,J)\) have opposite Bismut torsion three-forms, and their common torsion is closed:
\[
dd_I^c\omega_I=dd_J^c\omega_J=0.
\]
Thus, both Hermitian structures are pluriclosed (also called strong K\"ahler with torsion, or SKT). The ordinary K\"ahler case is recovered when \(I=J\) and \(\beta=0\). We adopt the conventions for \(d^c\) and the identification of differential forms with bundle maps used in \cite{ASU}. Thus, generalized Kähler geometry extends the classical Kähler framework by incorporating two complex structures, torsion, and holomorphic Poisson geometry, providing a natural setting in which to investigate canonical metrics and their rigidity properties.

Of particular interest is the class of generalized Kähler structures of symplectic type, which admits a description in terms of a symplectic form taming a complex structure; see \cite{GualtieriThesis,Gualtieri,ASU}. More precisely, fixing a complex manifold \((M,J)\), let \(F\) be a closed symplectic form satisfying \(F(X,JX)>0\) for every nonzero tangent vector \(X\). Such a form determines a Riemannian metric \(g\), a real two-form \(\beta\), and a second almost complex structure \(I\) by
\[
-FJ=g+\beta,\qquad I=-F^{-1}J^*F.
\]
If \(I\) is integrable, these data define a generalized Kähler structure of symplectic type \cite{ASU}. In particular, for a fixed complex structure \(J\), the entire bi-Hermitian structure is encoded by \(F\), subject to the additional integrability condition on \(I\).

An important feature of generalized Kähler geometry is its natural connection with Poisson geometry. As observed by Hitchin \cite{HitchinPoisson}, the commutator of the two complex structures determines a real Poisson bivector
\[
\pi=\frac12[I,J]g^{-1},
\]
whose \((2,0)\)-components with respect to \(I\) and \(J\) define holomorphic Poisson structures. In particular, the holomorphic Poisson tensor associated with $I$ and \(J\) are
\[
\pi_I=\pi-\sqrt{-1}I\pi, \quad \pi_J=\pi-\sqrt{-1}J\pi.
\]
This additional Poisson structure plays a central role in the study of generalized Kähler classes and canonical generalized Kähler metrics \cite{ASU}.

\subsection{Generalized K\"ahler classes}

In classical K\"ahler geometry, a central role is played by the space of K\"ahler metrics in a fixed K\"ahler class \(\alpha\in H^{1,1}(M,\mathbb R)\). By the \(\partial\bar\partial\)-lemma, fixing a reference K\"ahler form \(\omega_0\in\alpha\), every other K\"ahler form in the same class can be written as
\[
\omega_\varphi=\omega_0+dd_J^c\varphi,
\]
where \(\varphi\in C^\infty(M,\mathbb R)\) is unique up to an additive constant and satisfies \(\omega_\varphi>0\). Thus, a fixed K\"ahler class admits a global parametrization by an open convex space of smooth potentials. This provides the natural framework for the study of canonical K\"ahler metrics, including constant scalar curvature and extremal K\"ahler metrics.

The notion of a generalized K\"ahler class has emerged through the study of deformations of generalized K\"ahler structures, particularly the Hamiltonian deformation theory developed by Gualtieri \cite{GualtieriHamiltonian} and the deformation classes introduced by Gibson and Streets \cite{GibsonStreets}. Building on these developments, Apostolov, Streets, and Ustinovskiy \cite{ASU} formulated a notion of generalized K\"ahler class for structures of symplectic type with fixed complex structure \(J\), holomorphic Poisson tensor \(\pi_J\), and de Rham cohomology class \(\alpha\). They further established that the smooth path-connected components of the corresponding space of generalized K\"ahler structures can be characterized by Hamiltonian potential deformations, providing a natural framework for extending the Calabi program to generalized K\"ahler geometry.

More precisely, fix a compact complex manifold \((M,J)\), a holomorphic Poisson tensor \(\pi_J\), and a cohomology class \(\alpha\in H^2(M,\mathbb R)\). Denote by
\[
\mathcal{GK}_{\pi,\alpha}
=
\left\{
F\ \middle|
(F,J)\text{ is a generalized K\"ahler structure of symplectic type},
\ [F]=\alpha,\ \pi_J(F)=\pi_J
\right\}
\]
the space of generalized K\"ahler structures realizing these fixed data. Given a reference structure \(F_0\in\mathcal{GK}_{\pi,\alpha}\), the smooth path-connected component
\[
\mathcal{GK}_{\pi,\alpha}^{0}(F_0)
\]
containing \(F_0\) is called a \emph{generalized K\"ahler class} \cite{ASU}.

If we fix the complex structure \(J\), the holomorphic Poisson tensor \(\pi_J\), and the cohomology class \(\alpha\), infinitesimal variations within a generalized K\"ahler class can be represented by
\[
\dot F=dd_I^c\varphi,
\]
where \(\varphi\in C^\infty(M,\mathbb R)\) and \(I\) is the second complex structure associated with \(F\). Unlike the ordinary K\"ahler case, however, the complex structure \(I\) generally varies along the deformation. Consequently, the corresponding finite deformation cannot, in general, be expressed simply as \(F=F_0+dd_{I_0}^c\varphi\) for a single potential.

Instead, following the Hamiltonian deformation approach developed by Gualtieri \cite{GualtieriHamiltonian} and the formulation of Apostolov--Streets--Ustinovskiy \cite{ASU}, one considers a smooth family of real-valued functions \(\varphi_t\) and the associated evolution equations
\[
\frac{dF_t}{dt}=dd_{I_t}^c\varphi_t,
\qquad
I_t=(\Phi_t)_*I_0(\Phi_t)_*^{-1},
\]
where \(\Phi_t\) is the flow generated by the time-dependent vector field
\[
X_t=-\pi(d\varphi_t).
\]
Here \(\pi\) denotes the fixed real Poisson tensor determined by \(\pi_J\). Integrating the first equation gives
\[
F_t=F_0+\int_0^t dd_{I_s}^c\varphi_s\,ds.
\]
Since each variation is exact, the cohomology class remains unchanged:
\[
[F_t]=[F_0]=\alpha.
\]
Moreover, the deformation preserves the holomorphic Poisson tensor \(\pi_J\). Provided the taming condition remains satisfied, the family \(F_t\) therefore stays within the prescribed generalized K\"ahler class.

Apostolov--Streets--Ustinovskiy \cite{ASU} establishes that every generalized K\"ahler structure in the smooth path-connected component of \(F_0\) can be obtained through such a generalized K\"ahler potential deformation. In this sense, smooth real-valued functions play a role analogous to ordinary K\"ahler potentials, even though the deformation is nonlinear and involves the simultaneous evolution of \(F_t\) and \(I_t\). This provides a geometric characterization of generalized K\"ahler classes and makes them natural analogues of the classical spaces of K\"ahler potentials.

When the fixed holomorphic Poisson tensor vanishes, the real Poisson tensor \(\pi\) also vanishes, so the Hamiltonian vector field \(X_t\) is identically zero. Consequently, \(\Phi_t\) is the identity and \(I_t=I_0=J\). The deformation equations reduce to
\[
\frac{dF_t}{dt}=dd_J^c\varphi_t,
\]
and hence
\[
F_t=F_0+dd_J^c\left(\int_0^t\varphi_s\,ds\right).
\]
Thus, the generalized K\"ahler class reduces to the ordinary K\"ahler class,
\[
\boxed{
\mathcal{GK}_{0,\alpha}=\mathcal K_\alpha.
}
\]
In particular, generalized K\"ahler classes provide a natural extension of the classical framework for studying canonical metrics within fixed geometric data. Their potential-theoretic description also underlies the generalized Mabuchi functional and the variational formulation of the generalized Calabi problem developed in \cite{ASU}.

\subsection{Toric generalized K\"ahler classes and the two skew-symmetric parameters}

Let \((M^{2n},T)\) be a compact toric manifold, where \(T\) is an \(n\)-dimensional real torus with Lie algebra \(\mathfrak t\), and let \(P\subset\mathfrak t^*\cong\mathbb R^n\) be its Delzant moment polytope. We fix a collection of affine functions \(L=\{L_r\}_{r=1}^d\) defining the facets of \(P\), so that
\[
P=\{x\in\mathfrak t^*:L_r(x)\geq0,\quad r=1,\ldots,d\}.
\]
The primitive inward normals to these facets determine the integral lattice of the torus. Choosing an integral basis, we identify \(T\cong\mathbb R^n/(2\pi\mathbb Z)^n\). On the dense open orbit, we use momentum--angle coordinates
\[
(x,\theta)=(x_1,\ldots,x_n,\theta_1,\ldots,\theta_n)
\in P^\circ\times T.
\]

The Abreu--Guillemin description of toric generalized K\"ahler structures of symplectic type, developed in \cite{ASU} following earlier work on toric generalized K\"ahler geometry, associates to such a structure a triple
\[
(u,A,B),
\]
where \(u\) is a strictly convex symplectic potential on \(P^\circ\), defined up to an affine-linear function, and
\[
A,B\in\Lambda^2\mathfrak t
\]
are constant real skew-symmetric matrices. The ordinary toric K\"ahler geometry is recovered when \(A=B=0\). More generally, the potential \(u\) determines an associated ordinary toric K\"ahler structure, called the \emph{toric K\"ahler reduction}, while \(A\) and \(B\) encode two additional generalized K\"ahler deformations.

To describe these data explicitly, write
\[
G=\operatorname{Hess}u=(u_{ij}).
\]
In momentum--angle coordinates, the ordinary toric K\"ahler reduction is determined by the symplectic form
\[
\omega=\sum_{i=1}^n dx_i\wedge d\theta_i
\]
and the Riemannian metric
\[
g_u=\sum_{i,j=1}^n u_{ij}\,dx_i\,dx_j
+\sum_{i,j=1}^n u^{ij}\,d\theta_i\,d\theta_j,
\]
where \((u^{ij})=G^{-1}\). In particular, the corresponding Hessian metric on the interior of the moment polytope is \(G=\operatorname{Hess}u\).

The generalized K\"ahler structure is obtained by introducing the constant skew-symmetric matrices \(A\) and \(B\). In the conventions of \cite{ASU}, the complex structure \(J\) and the symplectic form \(F\) are expressed as
\begin{align}
J&=
\sum_{i,j=1}^n
\left[
(G+A)_{ij}\,dx_i\otimes
\frac{\partial}{\partial\theta_j}
-(G+A)^{-1}_{ij}\,d\theta_i\otimes
\frac{\partial}{\partial x_j}
\right],\label{eq:toric-J}\\
F&=
\sum_{i=1}^n dx_i\wedge d\theta_i
+\sum_{i,j=1}^n B_{ij}\,dx_i\wedge dx_j.
\label{eq:form}
\end{align}
The matrix \(G+A\) is invertible because its symmetric part \(G\) is positive definite. The skew-symmetric parameter \(A\) appears in the expression for \(J\), whereas \(B\) appears in the expression for \(F\). Notice that the sum in the second term of \eqref{eq:form} runs over all ordered pairs \((i,j)\), without a factor of \(1/2\). Thus, when \(n=2\) and
\[
B=\begin{pmatrix}0&t\\-t&0\end{pmatrix},
\]
we obtain
\[
\sum_{i,j=1}^2B_{ij}\,dx_i\wedge dx_j
=2t\,dx_1\wedge dx_2.
\]
This normalization will be important in the explicit examples considered later.

Although \(A\) and \(B\) enter the geometric construction in different ways, both contribute to the associated holomorphic Poisson tensor. More precisely, writing
\[
Z_i=\frac{\partial}{\partial\theta_i}
-\sqrt{-1}J\frac{\partial}{\partial\theta_i},
\]
Apostolov--Streets--Ustinovskiy \cite{ASU} shows
\[
\pi_J=\sum_{i,j=1}^n
(A_{ij}+\sqrt{-1}B_{ij})Z_i\wedge Z_j.
\]
Thus the holomorphic Poisson tensor is encoded by the complex skew-symmetric matrix \(A+\sqrt{-1}B\). For a fixed complex toric manifold and torus action, these coefficients uniquely determine the invariant holomorphic Poisson tensor. Consequently, fixing the holomorphic Poisson tensor within a toric generalized K\"ahler class fixes both \(A\) and \(B\).

Not every strictly convex symplectic potential is admissible for a prescribed \(B\). Let
\[
u_G(x)=\frac12\sum_{r=1}^d L_r(x)\log L_r(x)
\]
denote the Guillemin reference potential. According to \cite{ASU}, the admissible space \(\mathcal S_{A,B}(P,L)\) consists of potentials satisfying the Guillemin boundary condition
\[
u=u_G+v,\qquad v\in C^\infty(\overline P),
\]
together with the positivity condition
\[
\boxed{\operatorname{Hess}u+\sqrt{-1}B>0}
\]
on \(P^\circ\) and the corresponding tangential positivity conditions along every face of \(P\). Here \(\operatorname{Hess}u+\sqrt{-1}B\) is Hermitian because \(G\) is real symmetric and \(B\) is real skew-symmetric. The additional positivity requirement reflects the condition that the symplectic form \(F\) tames the complex structure \(J\), and ensures that the generalized K\"ahler structure extends smoothly to the compact toric manifold.

Finally, the generalized scalar curvature admits the Abreu-type expression
\begin{equation}
\label{eq:scalar}
\operatorname{Gscal}(u,A,B)
=
-\sum_{i,j=1}^n
\frac{\partial^2}{\partial x_i\partial x_j}
\left(\operatorname{Hess}u+\sqrt{-1}B\right)^{-1}_{ij}.
\end{equation}
In particular, the generalized scalar curvature equation depends on the potential \(u\) and the skew-symmetric matrix \(B\), but is independent of \(A\). This distinction plays a central role in our splitting results: when \(B\) respects a product decomposition of \(P\), the extremal equation forces the symplectic potential to split, even when \(A\) has nonzero mixed components. By contrast, mixed components of \(B\) can lead to extremal generalized K\"ahler structures with genuinely non-splitting symplectic potentials.

The admissible potentials form a convex set \(\Sc_B(P)\), independent of \(A\). The scalar equation is likewise independent of \(A\). These facts allow a theorem about \(u\) without a block condition on \(A\); they do not remove \(A\) from the underlying geometry.

\subsection{Results}

Our first result, Theorem \ref{thm:split}, establishes that every extremal toric generalized K\"ahler structure on a product toric manifold has a splitting symplectic potential whenever the skew-symmetric parameter \(B\) respects the product decomposition. Notably, no splitting assumption is imposed on the other skew-symmetric parameter \(A\).

We then construct two types of examples on \(\mathbb{CP}^1\times\mathbb{CP}^1\) demonstrating distinct failures of splitting. The first is an explicit constant-generalized-scalar-curvature structure with a nonzero mixed component of \(A\) and \(B=0\). Although its symplectic potential splits, the full generalized K"ahler structure does not. The second, established in Theorem \ref{thm:mixed}, consists of a one-parameter family of extremal generalized K\"ahler structures with nonzero mixed components of \(B\), whose symplectic potentials themselves fail to split. Moreover, Proposition \ref{prop:exactobstruction} provides an exact obstruction to the existence of splitting extremal potentials for nonzero mixed \(B\).

Together, these results distinguish the splitting of the symplectic potential from that of the full generalized K\"ahler structure and demonstrate the fundamentally different roles played by the two skew-symmetric parameters \(A\) and \(B\).

\begin{remark}
Our results are also related to the Yau--Tian--Donaldson type correspondence in toric generalized K\"ahler geometry. Apostolov, Streets, and Ustinovskiy \cite{ASU} proved that the existence of an extremal toric generalized K\"ahler structure implies the existence of an extremal K\"ahler metric in the same K\"ahler class. Conversely, they showed that an extremal toric K\"ahler metric gives rise to extremal generalized K\"ahler structures for sufficiently small multiples of any prescribed toric holomorphic Poisson tensor. Our splitting theorem complements these results by showing that, when \(B\) is block diagonal, an extremal generalized K\"ahler structure on a product toric manifold necessarily has a splitting symplectic potential, without any smallness assumption on the Poisson tensor. On the other hand, our examples with mixed \(B\) demonstrate that even an arbitrarily small Poisson deformation can produce extremal generalized K\"ahler structures with non-splitting symplectic potentials. Thus, while the existence of extremal K\"ahler and generalized K\"ahler structures is closely related, their splitting properties can differ substantially.
\end{remark}

\subsection*{Acknowledgment}
The first-named author would like to thank
Vestislav Apostolov for his interest in this work,
valuable suggestions, and enlightening discussions.

\subsection*{AI Usage Statement}

The authors used AI-assisted tools during the preparation
of this manuscript to assist with writing, editing, and
improving clarity. The authors have reviewed and verified
the content and take full responsibility for the accuracy,
originality, and integrity of the paper.

\section{Splitting of extremal potentials for block-diagonal \texorpdfstring{\(B\)}{B}}\label{sec:equation}

\subsection{Boundary conditions and conventions}

Let \(P\subset\mathfrak t^*\cong\mathbb R^n\) be a
compact Delzant polytope defined by affine functions
\(L_1,\ldots,L_d\), so that
\[
P=\{x\in\mathbb R^n:L_r(x)\geq 0,\ 
r=1,\ldots,d\}.
\]
We fix the facet labels and the corresponding
Lebesgue measure \(dx\) on \(\mathfrak t^*\).

Let \(B\in\Lambda^2\mathfrak t\) be a fixed real
skew-symmetric matrix. For a symplectic potential
\(u\), set
\[
H=\Hess u,\qquad
K_B(u)=H+ \sqrt{-1} B.
\]
Since \(H\) is real symmetric and \(B\) is real
skew-symmetric, \(K_B(u)\) is Hermitian.

Following \cite{ASU}, we denote by
\(\Sc_B(P)\) the space of admissible symplectic
potentials satisfying the following conditions:
\begin{enumerate}
\item \(u-u_G\in C^\infty(\overline P)\);
\item \(K_B(u)\) is positive definite on \(P^\circ\);
\item for every positive-dimensional face
\(Q\subset P\), the restriction
\[
\left.(\Hess u+\sqrt{-1} B)\right|_{TQ}
\]
is positive definite on the relative interior
\(Q^\circ\), where the restriction is understood
through the tangential Hessian of \(u|_Q\).
\end{enumerate}

An important feature of these conditions is that
they depend on \(B\), but not on the other
skew-symmetric parameter \(A\). In particular,
\cite{ASU} establishes
\[
\Sc_{A,B}(P,L)=\Sc_B(P),
\]
and shows that this space is convex. Thus,
if \(u_0,u_1\in\Sc_B(P)\), then
\[
u_s=(1-s)u_0+su_1\in\Sc_B(P),
\qquad 0\leq s\leq1.
\]

Define the real symmetric matrix
\[
X_B(u)=\Ree\bigl(K_B(u)^{-1}\bigr).
\]

For later use, we record an equivalent real
expression for \(X_B(u)\):
\begin{equation}\label{eq:real}
X_B(u)=
\bigl(H+BH^{-1}B\bigr)^{-1}.
\end{equation}
Indeed, writing
\(K_B(u)^{-1}=X+\sqrt{-1} Y\), with \(X\) symmetric
and \(Y\) skew-symmetric, the identity
\[
(H+\sqrt{-1} B)(X+\sqrt{-1} Y)=\mathrm{Id}
\]
gives
\[
HX-BY=\mathrm{Id},\qquad
HY+BX=0.
\]
Since \(H\) is invertible,
\(Y=-H^{-1}BX\), and hence
\[
(H+BH^{-1}B)X=\mathrm{Id},
\]
which proves \eqref{eq:real}. Notice that
\(BH^{-1}B\) is negative semidefinite, so the
positivity of \(K_B(u)\) imposes a nontrivial
restriction on the pair \((u,B)\).

We now formulate the extremal equation. Let
\(d\sigma\) denote the boundary measure determined
by the facet labels and \(dx\). The
\emph{extremal affine function}
\(\ell_P\in\Aff(P)\) is uniquely characterized by
\begin{equation}\label{eq:ell}
\int_P \ell_P a\,dx
=
2\int_{\partial P}a\,d\sigma,
\qquad
\text{for every }a\in\Aff(P).
\end{equation}
Existence and uniqueness follow from the
nondegeneracy of the \(L^2(P,dx)\) inner product
on the finite-dimensional space \(\Aff(P)\).

The generalized scalar curvature satisfies
the corresponding integration-by-parts identity,
so its \(L^2\)-projection onto \(\Aff(P)\) is
precisely \(\ell_P\), independently of the choice
of admissible potential.

A toric generalized K\"ahler structure is
extremal if and only if its generalized scalar
curvature is affine-linear on the moment
polytope. Thus, in the present notation, the
extremal equation is
\begin{equation}\label{eq:ext}
S_B(u)=\ell_P,
\qquad u\in\Sc_B(P).
\end{equation}
When \(\ell_P\) is constant, this reduces to
the constant generalized scalar curvature
equation.

Finally, we emphasize that both the admissible
space \(\Sc_B(P)\) and the extremal equation
\eqref{eq:ext} are independent of \(A\). This
observation will be important in the splitting
argument: the splitting behavior of the
symplectic potential is governed by \(B\), while
the parameter \(A\) may still contribute
nontrivially to the full generalized
K\"ahler geometry.

\subsection{Variational facts}
Now we define the relative Donaldson-Futaki invariant of $(P, L)$:
\[
\mathcal F_P(f)=2\int_{\partial P}f\,d\sigma-\int_P\ell_P f\,dx
\]
and, using any fixed reference \(u_*\in\Sc_0(P)\), the relative Mabuchi energy:
\begin{equation}\label{eq:energy}
\mathcal M_B(u)=\mathcal F_P(u)
-\int_P\log\det K_B(u)\,dx+\int_P\log\det\Hess u_*\,dx .
\end{equation}
All integrals are finite: the only singularities of the logarithms at facets are integrable logarithms. Reference changes add a constant.

For smooth \(f\) on \(\overline P\), the boundary conditions give
\begin{equation}\label{eq:ibp}
\int_P S_B(u)f\,dx
=2\int_{\partial P}f\,d\sigma
-\int_P\tr\bigl(X_B(u)\Hess f\bigr)\,dx.
\end{equation}
Differentiating \eqref{eq:energy} and applying \eqref{eq:ibp} yields
\[
D\mathcal M_B(u)[f]=\int_P(S_B(u)-\ell_P)f\,dx .
\]
For a segment \(u_s=u+sf\) of admissible potentials,
\begin{equation}\label{eq:convex}
\frac{d^2}{ds^2}\mathcal M_B(u_s)
=\int_P\tr\left[
 \left(K_B(u_s)^{-1/2}(\Hess f)K_B(u_s)^{-1/2}\right)^2
 \right]dx\geq0.
\end{equation}
The matrix inside the square is Hermitian; hence the trace of its square is a sum of squares of real eigenvalues. Equality implies \(\Hess f=0\). The boundary estimates implicit in the admissibility conditions justify differentiation and integration.

It follows directly that any solution of \eqref{eq:ext} globally minimizes \(\mathcal M_B\) on the convex admissible set, and that any two solutions differ by an affine function. This is a uniqueness statement for fixed toric coordinates and fixed \(B\).

\subsection{The product averaging argument}
\label{sec:split}

In this subsection, we establish our main splitting result
for extremal toric generalized K\"ahler potentials.
The proof is based on an averaging procedure over the
factors of a product polytope and the convexity properties
of the generalized Mabuchi functional.

Let \(P_1\subset\mathbb R^{n_1}\) and
\(P_2\subset\mathbb R^{n_2}\) be compact Delzant polytopes,
and consider their product
\[
P=P_1\times P_2
\subset\mathbb R^{n_1}\times\mathbb R^{n_2}.
\]
We equip \(P\) with the product facet labels and write
\((x,y)\) for its coordinates, where \(x\in P_1\) and
\(y\in P_2\). Denote the volumes of the factors by
\[
V_1=\int_{P_1}dx,\qquad
V_2=\int_{P_2}dy.
\]

We assume throughout this subsection that the constant
skew-symmetric matrix \(B\) respects the product
decomposition, namely,
\[
B=
\begin{pmatrix}
B_1&0\\
0&B_2
\end{pmatrix},
\qquad
B_i\in\Lambda^2\mathbb R^{n_i}.
\]
No corresponding assumption will be imposed on the
other skew-symmetric parameter \(A\).

We first observe that the extremal affine function
is compatible with products.

\begin{lemma}\label{lem:ell}
Let \(P=P_1\times P_2\) be a product of compact
Delzant polytopes equipped with the product facet
labels. Then
\[
\ell_P(x,y)=\ell_{P_1}(x)+\ell_{P_2}(y).
\]
\end{lemma}

\begin{proof}
Recall that the extremal affine function on a
labeled polytope \(Q\) is uniquely characterized by
\[
\int_Q \ell_Q a\,dx
=
2\int_{\partial Q}a\,d\sigma_Q,
\qquad a\in\Aff(Q).
\]

Set
\[
V_i=\operatorname{Vol}(P_i),
\qquad
S_i=\int_{\partial P_i}d\sigma_i.
\]
In particular,
\[
\int_{P_i}\ell_{P_i}\,dx_i=2S_i.
\]

We claim that
\[
\ell(x,y)=\ell_{P_1}(x)+\ell_{P_2}(y)
\]
satisfies the defining identity for \(\ell_P\).

Every affine function on \(P\) can be expressed as
\[
a(x,y)=a_1(x)+a_2(y),
\]
where \(a_i\in\Aff(P_i)\). The product structure
of the volume and boundary measures gives
\[
\begin{aligned}
\int_P\ell a
={}&
V_2\int_{P_1}\ell_{P_1}a_1
+V_1\int_{P_2}\ell_{P_2}a_2\\
&+
2S_2\int_{P_1}a_1
+2S_1\int_{P_2}a_2.
\end{aligned}
\]
Using the defining identities for
\(\ell_{P_1}\) and \(\ell_{P_2}\), this becomes
\[
\begin{aligned}
\int_P\ell a
={}&
2V_2\int_{\partial P_1}a_1\,d\sigma_1
+2V_1\int_{\partial P_2}a_2\,d\sigma_2\\
&+
2S_2\int_{P_1}a_1
+2S_1\int_{P_2}a_2.
\end{aligned}
\]

On the other hand, the facets of \(P\) are
precisely \(E_1\times P_2\) and
\(P_1\times E_2\), where \(E_i\) ranges over
the facets of \(P_i\). Since the facet measures
are the corresponding product measures,
\[
\begin{aligned}
2\int_{\partial P}a\,d\sigma
={}&
2V_2\int_{\partial P_1}a_1\,d\sigma_1
+2V_1\int_{\partial P_2}a_2\,d\sigma_2\\
&+
2S_2\int_{P_1}a_1
+2S_1\int_{P_2}a_2.
\end{aligned}
\]

Consequently,
\[
\int_P\ell a
=
2\int_{\partial P}a\,d\sigma
\]
for every \(a\in\Aff(P)\). By uniqueness of
the extremal affine function,
\[
\ell_P=\ell_{P_1}+\ell_{P_2}.
\]
\end{proof}

We next introduce an averaging operator on the
space of symplectic potentials.

For \(u\in\Sc_B(P)\), define its partial averages by
\[
m_1(x)=\frac1{V_2}\int_{P_2}u(x,y)\,dy,
\qquad
m_2(y)=\frac1{V_1}\int_{P_1}u(x,y)\,dx,
\]
and set
\[
c_u=\frac1{V_1V_2}\int_Pu(x,y)\,dx\,dy.
\]
The \emph{product averaging operator} is defined by
\begin{equation}\label{eq:projection}
(Tu)(x,y)=m_1(x)+m_2(y)-c_u.
\end{equation}
The normalization is chosen so that
\[
\frac1{V_2}\int_{P_2}(Tu)(x,y)\,dy=m_1(x)
\]
and
\[
\frac1{V_1}\int_{P_1}(Tu)(x,y)\,dx=m_2(y).
\]
In particular, \(T\) preserves the partial averages
of \(u\), and
\[
T^2u=Tu.
\]
Thus \(T\) is a projection onto the space of
potentials that split as sums of functions on
the two factors.

Our first task is to verify that this projection
preserves generalized K\"ahler admissibility.

\begin{lemma}[Admissibility of averaging]
\label{lem:admiss}
Suppose \(B=B_1\oplus B_2\). If
\(u\in\Sc_B(P)\), then
\[
m_1\in\Sc_{B_1}(P_1),
\qquad
m_2\in\Sc_{B_2}(P_2),
\]
and consequently
\[
Tu\in\Sc_B(P).
\]
\end{lemma}

\begin{proof}
We verify the boundary regularity and positivity
conditions separately.

Since \(P=P_1\times P_2\), its Guillemin potential
splits:
\[
u_G(x,y)=u_{G,1}(x)+u_{G,2}(y).
\]
By admissibility,
\[
u(x,y)=u_{G,1}(x)+u_{G,2}(y)+h(x,y),
\qquad
h\in C^\infty(\overline P).
\]
Averaging over \(P_2\), we obtain
\[
m_1(x)=u_{G,1}(x)+c_2+
\frac1{V_2}\int_{P_2}h(x,y)\,dy,
\]
where \(c_2\) is a constant. Hence
\[
m_1-u_{G,1}\in C^\infty(\overline P_1).
\]
The same argument applies to \(m_2\).

We next consider positivity in the interior.
Write
\[
\Hess u+\sqrt{-1} B=
\begin{pmatrix}
H_{xx}+\sqrt{-1} B_1&H_{xy}\\
H_{yx}&H_{yy}+\sqrt{-1}B_2
\end{pmatrix}.
\]
Since this matrix is positive definite, its
principal diagonal blocks are positive definite.
Moreover,
\[
\Hess m_1+\sqrt{-1} B_1
=
\frac1{V_2}\int_{P_2}
\bigl(H_{xx}+\sqrt{-1} B_1\bigr)\,dy.
\]
The right-hand side is an average of positive
definite Hermitian matrices, and is therefore
positive definite. Similarly,
\[
\Hess m_2+\sqrt{-1} B_2>0.
\]

It remains to verify positivity along the faces.
Let \(E_1\subset P_1\) be a positive-dimensional
face. On the relative interior of
\(E_1\times P_2\), the admissibility of \(u\)
implies that the tangential restriction of
\(\Hess u+\sqrt{-1} B\) is positive definite.
Restricting further to vectors tangent to
\(E_1\), and averaging over \(P_2\), yields
\[
\left.(\Hess m_1+\sqrt{-1} B_1)\right|_{TE_1}>0.
\]
The tangential Hessians are well defined because
the restrictions of the smooth parts of the
potentials extend to the faces. The corresponding
argument applies to every face of \(P_2\).

Thus \(m_i\in\Sc_{B_i}(P_i)\). Finally,
\[
\Hess(Tu)+\sqrt{-1}B=
\begin{pmatrix}
\Hess m_1+\sqrt{-1} B_1&0\\
0&\Hess m_2+\sqrt{-1} B_2
\end{pmatrix},
\]
which satisfies all the admissibility conditions
on the product polytope.
\end{proof}

The next lemma is the central variational
ingredient. It shows that averaging decreases
the generalized Mabuchi functional, with strict
inequality unless the potential already splits.

\begin{lemma}[Energy decreases under averaging]
\label{lem:decrease}
Suppose \(B=B_1\oplus B_2\). Then, for every
\(u\in\Sc_B(P)\),
\[
\mathcal M_B(Tu)\leq\mathcal M_B(u).
\]
Moreover, equality holds if and only if
\[
\Hess_{xy}u=0
\]
throughout the interior of \(P\).
\end{lemma}

\begin{proof}
Recall from \eqref{eq:energy} that
\[
\mathcal M_B(u)
=
\mathcal F_P(u)
-\int_P\log\det(\Hess u+\sqrt{-1} B)\,dx\,dy
+C,
\]
where \(C\) is independent of \(u\).

We first compare the linear terms.
The equality of partial averages implies
\[
\int_PTu\,dx\,dy=\int_Pu\,dx\,dy.
\]
More generally, Lemma \ref{lem:ell} and the
partial-average identities give
\[
\int_P\ell_P(Tu)\,dx\,dy
=
\int_P\ell_Pu\,dx\,dy.
\]
Furthermore, the boundary of \(P\) is the union
of the product facets
\[
\partial P=
(\partial P_1\times P_2)
\cup(P_1\times\partial P_2),
\]
up to sets of boundary measure zero.
Thus,
\[
\int_{\partial P}Tu\,d\sigma
=
\int_{\partial P}u\,d\sigma.
\]
Consequently,
\begin{equation}\label{eq:linear-preserved}
\mathcal F_P(Tu)=\mathcal F_P(u).
\end{equation}

It remains to compare the determinant terms.
Set
\[
K=\Hess u+\sqrt{-1}B
=
\begin{pmatrix}
D_1&C\\
C^*&D_2
\end{pmatrix},
\]
where
\[
D_1=H_{xx}+\sqrt{-1}B_1,\qquad
D_2=H_{yy}+\sqrt{-1}B_2,\qquad
C=H_{xy}.
\]
Since \(K\) is positive definite Hermitian,
so are \(D_1\) and \(D_2\).

By the Schur complement formula,
\[
\det K=
\det D_1\,
\det(D_2-C^*D_1^{-1}C).
\]
Since \(C^*D_1^{-1}C\) is positive semidefinite,
the Fischer determinant inequality yields
\[
\det K\leq\det D_1\det D_2.
\]
Equality holds precisely when \(C=0\).
Taking logarithms, we obtain
\begin{equation}\label{eq:fischer}
\log\det K
\leq\log\det D_1+\log\det D_2.
\end{equation}

We next apply the concavity of \(\log\det\)
on the cone of positive definite Hermitian
matrices. For each fixed \(x\in P_1^\circ\),
Jensen's inequality gives
\[
\begin{aligned}
\frac1{V_2}\int_{P_2}\log\det D_1(x,y)\,dy
&\leq
\log\det\left(
\frac1{V_2}\int_{P_2}D_1(x,y)\,dy
\right)\\
&=
\log\det(\Hess m_1(x)+\sqrt{-1}B_1).
\end{aligned}
\]
Integrating over \(P_1\), we conclude that
\[
\int_P\log\det D_1\,dx\,dy
\leq
V_2\int_{P_1}
\log\det(\Hess m_1+\sqrt{-1}B_1)\,dx.
\]
Similarly,
\[
\int_P\log\det D_2\,dx\,dy
\leq
V_1\int_{P_2}
\log\det(\Hess m_2+\sqrt{-1}B_2)\,dy.
\]

Combining these inequalities with
\eqref{eq:fischer}, we obtain
\[
\begin{aligned}
\int_P\log\det K\,dx\,dy
&\leq
V_2\int_{P_1}
\log\det(\Hess m_1+\sqrt{-1}B_1)\,dx\\
&\quad+
V_1\int_{P_2}
\log\det(\Hess m_2+\sqrt{-1}B_2)\,dy\\
&=
\int_P\log\det(\Hess(Tu)+\sqrt{-1}B)\,dx\,dy.
\end{aligned}
\]
The logarithmic boundary singularities are
integrable, so these inequalities extend to
the integrals over the full polytopes.

Together with \eqref{eq:linear-preserved},
this proves
\[
\mathcal M_B(Tu)\leq\mathcal M_B(u).
\]

Finally, if equality holds, the nonnegative
difference in the pointwise Fischer inequality
must vanish almost everywhere. Hence
\(C=H_{xy}=0\) almost everywhere, and therefore
everywhere in \(P^\circ\) by smoothness.
Conversely, if \(H_{xy}=0\), then \(u\) is a sum
of functions on the two factors, up to a
constant, and \(Tu=u\). Equality follows.
\end{proof}

We are now ready to prove the splitting theorem.

\begin{theorem}[Splitting of extremal potentials]
\label{thm:split}
Let \(P=P_1\times P_2\) be a product of compact
Delzant polytopes, and suppose that
\(B=B_1\oplus B_2\).
If \(u\in\Sc_B(P)\) is an extremal toric
generalized K\"ahler potential, that is,
\[
S_B(u)=\ell_P,
\]
then \(u\) splits, up to an affine-linear
function:
\[
u(x,y)=u_1(x)+u_2(y)+a(x,y),
\qquad a\in\Aff(P),
\]
where
\[
u_i\in\Sc_{B_i}(P_i),
\qquad
S_{B_i}(u_i)=\ell_{P_i}.
\]
No block-diagonal assumption on \(A\) is required.
\end{theorem}

\begin{proof}
Since \(u\) solves the extremal equation, it is
a critical point of the relative generalized
Mabuchi functional. By the convexity established
in \eqref{eq:convex}, it is a global minimizer
on the admissible space \(\Sc_B(P)\).

By Lemma \ref{lem:admiss}, the averaged potential
\(Tu\) is also admissible. Therefore,
\[
\mathcal M_B(u)\leq\mathcal M_B(Tu).
\]
On the other hand, Lemma \ref{lem:decrease} gives
\[
\mathcal M_B(Tu)\leq\mathcal M_B(u).
\]
Hence equality holds. The equality
characterization in Lemma \ref{lem:decrease}
implies
\[
\Hess_{xy}u=0.
\]

Since \(P^\circ=P_1^\circ\times P_2^\circ\)
is connected, this identity implies that
\(u\) is the sum of functions depending only
on the respective factors. In particular,
\[
u=Tu=m_1(x)+m_2(y)-c_u.
\]
The admissibility of the factors follows from
Lemma \ref{lem:admiss}.

We next verify the extremal equations on the
factors. Since both \(u\) and \(B\) split,
\[
\Hess u+\sqrt{-1}B=
\begin{pmatrix}
\Hess m_1+\sqrt{-1}B_1&0\\
0&\Hess m_2+\sqrt{-1}B_2
\end{pmatrix}.
\]
The generalized Abreu formula therefore gives
\[
S_B(u)=S_{B_1}(m_1)+S_{B_2}(m_2).
\]
Together with Lemma \ref{lem:ell}, this yields
\[
S_{B_1}(m_1)(x)-\ell_{P_1}(x)
=
-\bigl(S_{B_2}(m_2)(y)-\ell_{P_2}(y)\bigr).
\]
Since the two sides depend on independent
variables, each must be constant.

The integration-by-parts identity \eqref{eq:ibp},
applied to the constant test function on each
factor, shows that these constants vanish.
Thus
\[
S_{B_i}(m_i)=\ell_{P_i},
\qquad i=1,2.
\]

Finally, the admissible space and the generalized
scalar curvature equation are independent of
\(A\). Hence no splitting hypothesis on \(A\)
is needed.
\end{proof}

\begin{corollary}
Suppose \(B=B_1\oplus B_2\). Then an extremal
toric generalized K\"ahler potential exists
on \(P_1\times P_2\) if and only if extremal
potentials exist on both factors.
\end{corollary}

\begin{proof}
The forward implication follows from
Theorem \ref{thm:split}. Conversely, if
\(u_i\in\Sc_{B_i}(P_i)\) solve the corresponding
extremal equations, then
\[
u(x,y)=u_1(x)+u_2(y)
\]
belongs to \(\Sc_B(P)\), and the additivity
of generalized scalar curvature and
Lemma \ref{lem:ell} show that it is extremal.
\end{proof}

\begin{remark}[Splitting of the geometric structure]
The conclusion of Theorem \ref{thm:split}
concerns the symplectic potential and,
consequently, the associated ordinary toric
K\"ahler reduction. It does not necessarily
imply splitting of the full generalized
K\"ahler structure. Indeed, even when
\(\Hess u\) and \(B\) are block diagonal,
mixed components of \(A\) may couple the
two factors through the matrix \(\Hess u+A\).
If both \(A\) and \(B\) are block diagonal,
the full toric generalized K\"ahler structure
splits. The theorem extends inductively
to products of finitely many factors.
\end{remark}

\section{An explicit non-product example with split potential}\label{sec:gkproduct}

Potential splitting and splitting of the full generalized K\"ahler structure
are distinct questions. We now give an intrinsic obstruction to the latter.

\subsection{Meaning of a GK product}
A product of GK structures on two curves has bi-Hermitian tensors
\[
g=g_1\oplus g_2,\qquad I=I_1\oplus I_2,\qquad
J=J_1\oplus J_2,\qquad
\beta=\operatorname{pr}_1^*\beta_1+\operatorname{pr}_2^*\beta_2.
\]
A GK product decomposition up to diffeomorphism means that a diffeomorphism
identifies these tensors with the given structure. The obstruction below
also survives closed two-form transforms, since these change \(\beta\)
but not \(g,I,J\). We do not claim here that the Riemannian metric alone
has no product decomposition.

\begin{lemma}\label{lem:curvepoisson}
Every GK structure on a real surface has
\([I,J]=0\), and hence has zero canonical real Poisson tensor
\(\pi=\tfrac12[I,J]g^{-1}\). A product of two GK curves therefore has
\(\pi=0\) and \(\pi_J=0\).
\end{lemma}
\begin{proof}
On a two-dimensional Euclidean vector space the orthogonal complex
structures are exactly the two rotations by \(+\pi/2\) and \(-\pi/2\).
Thus two complex structures Hermitian for the same metric satisfy
\(I=J\) or \(I=-J\) at every point, and commute in either case.
On a product the commutator is the direct sum of the two factor
commutators. The assertion follows, as does its invariance under
diffeomorphism.
\end{proof}

\subsection{A concrete cscGK non-product with split potential}
There is also a fully explicit example with a fixed numerical skew
parameter and no perturbative existence argument. It uses the
\emph{real skew addition to the complex-structure matrix}, denoted
by \(A\), rather than the parameter \(B\) in \(H+\sqrt{-1}B\).

Let us consider $X = \mathbb{CP}^1 \times \mathbb{CP}^1$, the corresponding Delzant polytope is the standard square.
We let
\[
w(x,y)=f(x)+f(y),\qquad
f(s)=\tfrac12\bigl((1+s)\log(1+s)+(1-s)\log(1-s)\bigr),
\]
and take
\begin{equation}\label{eq:explicitA}
A=\begin{pmatrix}0&1/10\\-1/10&0\end{pmatrix},\qquad B=0.
\end{equation}
These data are globally admissible: \(w\) has the standard Guillemin
boundary behavior, and admissibility is independent of \(A\). The matrix in the toric complex-structure
ansatz is explicitly
\[
D^2w+A=
\begin{pmatrix}
(1-x^2)^{-1}&1/10\\
-1/10&(1-y^2)^{-1}
\end{pmatrix}.
\]
By \eqref{eq:scalar}, its generalized scalar curvature is
\[
S_0(w)=-\partial_x^2(1-x^2)-\partial_y^2(1-y^2)=4.
\]
On the other hand, let \(X_i=\partial/\partial\theta_i\), \(i=1,2\),
denote the infinitesimal generators of the real torus
action on \(\mathbb{CP}^1\times\mathbb{CP}^1\).
Define
\[
Z_i=X_i-\sqrt{-1}  JX_i,\qquad i=1,2.
\]
These are holomorphic vector fields of type \((1,0)\)
associated with the complexified torus action.
On the dense open orbit \((\mathbb C^*)^2\), in
holomorphic coordinates \((z_1,z_2)\) adapted to
the torus action, they take the form
\[
Z_i=2 \sqrt{-1} z_i\frac{\partial}{\partial z_i}.
\]

Thus the corresponding holomorphic Poisson tensor is
\[
\pi_J
=\sum_{i,j=1}^2 A_{ij}Z_i\wedge Z_j
=\frac{1}{5} \,Z_1\wedge Z_2.
\]
Equivalently, in holomorphic coordinates,
\[
\pi_J
=-\frac{4}{5}\,z_1z_2
\frac{\partial}{\partial z_1}
\wedge
\frac{\partial}{\partial z_2}.
\]
This tensor has a nonzero mixed
component with respect to the product decomposition.
Consequently, the generalized K\"ahler structure
does not split as a product, even though its
symplectic potential is a sum of potentials
on the two factors.

\section{Mixed-\texorpdfstring{\(B\)}{B} deformations with non-split extremal potentials}\label{sec:deformation}

In the previous section, we constructed an extremal toric generalized K\"ahler structure on \(\mathbb{CP}^1\times\mathbb{CP}^1\) with a splitting symplectic potential but a non-splitting generalized K\"ahler structure, arising from a nonzero mixed component of the skew-symmetric parameter \(A\). This shows that even when the symplectic potential splits, mixed components of \(A\) can obstruct splitting of the full generalized K\"ahler structure. To investigate whether the symplectic potential itself must split, we now eliminate this obstruction by assuming that \(A\) respects the product decomposition. In fact, we take \(A=0\) throughout this section and instead introduce a nonzero mixed component of \(B\). We show that, even under this restriction on \(A\), there exist extremal toric generalized K\"ahler structures on \(\mathbb{CP}^1\times\mathbb{CP}^1\) whose symplectic potentials do not split.

\subsection{Local deformation from a K\"ahler extremal metric}
\label{sec:ift}

We begin by recalling the local deformation theory
for extremal toric generalized K\"ahler structures
developed by Apostolov--Streets--Ustinovskiy
\cite{ASU}. Their argument adapts the
openness method of LeBrun--Simanca \cite{LS}
to the generalized Abreu equation. The key
observation is that, near an ordinary extremal
K\"ahler metric, the generalized extremal
equation can be solved by the implicit function
theorem after removing the finite-dimensional
kernel arising from holomorphic vector fields.

Let \(P\subset\mathbb R^n\) be a compact Delzant
polytope with fixed facet labels, and let
\(u_0\in\Sc_0(P)\) be the symplectic potential
of a smooth toric extremal K\"ahler metric.
Thus,
\[
S_0(u_0)=\ell_P,
\]
where \(\ell_P\) is the extremal affine function
determined by the labeled polytope.

Fix a real skew-symmetric matrix
\(B_*\in\Lambda^2\mathbb R^n\), and consider
the one-parameter family
\[
B_t=tB_*.
\]
We seek admissible symplectic potentials
\(u_t\in\Sc_{B_t}(P)\) satisfying
\begin{equation}\label{eq:deformed-extremal}
S_{B_t}(u_t)=\ell_P,
\qquad u_{t=0}=u_0.
\end{equation}
The generalized scalar curvature is given by
\[
S_{B_t}(u)
=
-\sum_{i,j=1}^n
\frac{\partial^2}{\partial x_i\partial x_j}
\left(\Hess u+\sqrt{-1}tB_*\right)^{-1}_{ij}.
\]
The equation depends on \(B_t\), but not on
the other skew-symmetric parameter \(A\).

\begin{proposition}[Local deformation of extremal potentials]
\label{prop:ift}
Let \(u_0\in\Sc_0(P)\) be the symplectic
potential of a smooth compact toric extremal
K\"ahler metric, and let
\(B_*\in\Lambda^2\mathbb R^n\) be fixed.
Then there exists \(\varepsilon>0\) and a
smooth family of admissible potentials
\[
u_t\in\Sc_{tB_*}(P),
\qquad |t|<\varepsilon,
\]
such that
\[
S_{tB_*}(u_t)=\ell_P,
\qquad u_{t=0}=u_0.
\]
The potentials may be normalized by requiring
\[
u_t-u_0\perp\Aff(P)
\]
with respect to the \(L^2(P,dx)\) inner product.
With this normalization, the family is locally
unique and satisfies
\[
u_{-t}=u_t.
\]
Consequently, it admits an expansion
\begin{equation}\label{eq:even-expansion}
u_t=u_0+t^2v+O(t^4),
\end{equation}
where the expansion holds in every fixed
finite smooth norm on the compact toric
manifold, using the corresponding
K\"ahler-potential identification.
\end{proposition}

\begin{proof}
The local existence statement follows from
\cite{ASU}. We recall the
implicit-function-theorem construction and
explain why the normalized family can be
chosen to depend smoothly on \(t\).

Let \((M,J,\omega_0)\) be the compact toric
K\"ahler manifold associated with \(u_0\).
Instead of working directly with symplectic
potentials, which have logarithmic
singularities along the boundary of \(P\),
we use the corresponding smooth
\(T\)-invariant K\"ahler potentials on \(M\).

For a sufficiently small invariant function
\(\varphi\), set
\[
\omega_\varphi
=\omega_0+dd_J^c\varphi,
\]
and let \(\mu^\varphi\) be the corresponding
moment map, normalized to have image \(P\).
Following \cite{ASU}, define
\[
F_{\varphi,t}
=
\omega_\varphi+
t\sum_{i,j=1}^n(B_*)_{ij}\,
d\mu_i^\varphi\wedge d\mu_j^\varphi.
\]
For \(\varphi\) and \(t\) sufficiently small,
this form tames \(J\) and determines a toric
generalized K\"ahler structure of symplectic
type. Under the Legendre-transform
identification, its Abreu--Guillemin data
are \((u,0,tB_*)\).

Consider the nonlinear operator
\[
W(\varphi,t)
=
\operatorname{Gscal}(F_{\varphi,t},J)
-\ell_P(\mu^\varphi).
\]
The generalized extremal equation is
equivalent to
\[
W(\varphi,t)=0.
\]
At the reference metric,
\[
W(0,0)=0.
\]

We first establish the regularity of this
operator. Fix \(k\geq0\) and
\(0<\gamma<1\). On a sufficiently small
open neighborhood of \((0,0)\), the
operator \(W\) extends to a \(C^\infty\) map
\[
W:
C^{k+4,\gamma}(M,\mathbb R)^T
\times(-\varepsilon_0,\varepsilon_0)
\longrightarrow
C^{k,\gamma}(M,\mathbb R)^T.
\]

At \((0,0)\), the linearization of \(W\)
in the potential direction is the
ordinary Lichnerowicz operator:
\[
D_\varphi W(0,0)=\mathcal L_{\omega_0}.
\]
Since \(\omega_0\) is extremal, this
operator is elliptic and self-adjoint.
Its kernel on \(T\)-invariant functions
consists of the affine Hamiltonians,
namely,
\[
\mathcal K
=
\{a\circ\mu^0:a\in\Aff(P)\}.
\]

Let \(\Pi_0\) denote the
\(L^2(M,\omega_0^n/n!)\)-orthogonal
projection onto \(\mathcal K\).
For each \(k\geq0\), define
\[
\begin{aligned}
E_k&=
C^{k+4,\gamma}(M,\mathbb R)^T
\cap\mathcal K^\perp,\\
Y_k&=
C^{k,\gamma}(M,\mathbb R)^T
\cap\mathcal K^\perp.
\end{aligned}
\]
Consider the projected operator
\[
\mathcal F(\varphi,t)
=
(\mathrm{Id}-\Pi_0)W(\varphi,t).
\]
It defines a \(C^\infty\) map
\[
\mathcal F:E_k\times
(-\varepsilon_0,\varepsilon_0)
\longrightarrow Y_k.
\]
Moreover,
\[
D_\varphi\mathcal F(0,0)
=
\mathcal L_{\omega_0}|_{E_k}.
\]
By elliptic Fredholm theory and
Schauder estimates,
\[
\mathcal L_{\omega_0}:E_k\longrightarrow Y_k
\]
is a bounded linear isomorphism.

The smooth implicit function theorem
therefore provides a unique sufficiently
small smooth curve
\[
t\longmapsto\varphi_t\in E_k,
\qquad \varphi_0=0,
\]
satisfying
\[
(\mathrm{Id}-\Pi_0)W(\varphi_t,t)=0.
\]

We next verify that the projected solution
solves the full extremal equation.
The integration-by-parts characterization
of the extremal affine function implies
\[
\Pi_tW(\varphi_t,t)=0,
\]
where \(\Pi_t\) is the orthogonal projection
onto the affine Hamiltonians associated
with the deformed moment map
\(\mu^{\varphi_t}\), using the corresponding
symplectic volume form.

Since the projected equation implies
\[
W(\varphi_t,t)\in\mathcal K,
\]
and the pairing between the fixed affine
Hamiltonians and their deformed
counterparts remains nondegenerate for
small \(t\), we conclude that
\[
W(\varphi_t,t)=0.
\]

Thus \(\varphi_t\) defines a family of
extremal generalized K\"ahler structures.
Elliptic regularity shows that each
\(\varphi_t\) is smooth on \(M\).

Furthermore, the implicit function theorem
can be applied in every higher
H\"older space. The resulting local
solutions agree by uniqueness in the
normalized slice. Consequently, the
family depends smoothly on \(t\) in
all \(C^k\)-norms.

The inverse Legendre transform gives
the corresponding symplectic potentials
\(u_t\). Since the Legendre-transform
identification depends smoothly on
the K\"ahler potential near the reference
metric, we obtain
\[
t\longmapsto u_t-u_0
\in C^\infty(\overline P)
\]
smoothly. The potentials are defined
up to affine-linear functions, which
we fix by requiring
\[
u_t-u_0\perp\Aff(P).
\]
This normalization preserves smooth
parameter dependence.

Finally, we establish the symmetry
under \(t\mapsto-t\). For every admissible
potential \(u\), we have
\[
\Hess u- tB_* \sqrt{-1}
=
\overline{\Hess u+ tB_* \sqrt{-1}},
\]
and hence
\[
S_{-tB_*}(u)=S_{tB_*}(u).
\]
The admissibility conditions are
also invariant under \(t\mapsto-t\).

Therefore, the normalized potentials
\(u_t\) and \(u_{-t}\) satisfy the same
generalized extremal equation for
each sufficiently small \(t\).
Local uniqueness implies
\[
u_{-t}=u_t.
\]

The normalized family is thus smooth
and even in \(t\). In particular,
all its odd derivatives at \(t=0\)
vanish. Taylor's theorem gives
\[
u_t=u_0+t^2v+O_{C^k}(t^4),
\qquad
v=\frac12
\left.\frac{d^2u_t}{dt^2}\right|_{t=0},
\]
for every fixed \(k\geq0\).
\end{proof}

\subsection{Non-splitting extremal potentials on
\texorpdfstring{$\mathbb{CP}^1\times\mathbb{CP}^1$}
{CP1 x CP1}}
\label{sec:mixedexample}

We now apply the deformation theory of the previous
subsection to construct a family of extremal toric
generalized K\"ahler structures on
\(\mathbb{CP}^1\times\mathbb{CP}^1\) whose symplectic
potentials do not split. Starting from the product
Fubini--Study metric, we introduce a mixed
skew-symmetric parameter \(B_t\). We first determine
the leading correction to the symplectic potential
and then establish an exact obstruction to splitting.

To investigate whether the extremal potentials
constructed in Proposition~\ref{prop:ift} split,
we examine their leading-order deformation
from the product Fubini--Study potential \(u_0\).

Let
\[
\mathcal S(u,t)=S_{tB_*}(u),
\qquad
B_*=
\begin{pmatrix}
0&1\\
-1&0
\end{pmatrix}.
\]
By Proposition~\ref{prop:ift}, the normalized
family \(u_t\) is smooth and even in \(t\).
Consequently, Taylor's theorem gives
\[
u_t=u_0+t^2v+O(t^4),
\qquad
v=\frac12
\left.\frac{d^2u_t}{dt^2}\right|_{t=0},
\]
where the expansion is understood in the
appropriate smooth norms for \(u_t-u_0\).

For the fixed potential \(u_0\), the generalized
scalar curvature is an even function of \(t\).
Since \(S_0(u_0)=4\), we may write
\[
\mathcal S(u_0,t)
=4-t^2Q+O(t^4),
\]
where
\[
Q=-\frac12
\left.\frac{\partial^2}{\partial t^2}
\mathcal S(u_0,t)\right|_{t=0}.
\]
The function \(Q\) measures the leading-order
change in generalized scalar curvature caused
by the deformation of \(B\), while keeping
the symplectic potential fixed.

We now account for the variation of the potential.
Taylor expansion of \(\mathcal S\) with respect
to \(u\), around \(u_0\), gives
\[
\begin{aligned}
\mathcal S(u_t,t)
&=\mathcal S(u_0,t)
+D_u\mathcal S(u_0,t)[u_t-u_0]\\
&\quad+O\bigl(\|u_t-u_0\|^2\bigr).
\end{aligned}
\]
At \(t=0\), the derivative with respect to \(u\)
is the linearization of the ordinary Abreu
scalar curvature operator:
\[
D_u\mathcal S(u_0,0)[v]=Lv.
\]
Since \(\mathcal S(u,t)\) is even in \(t\)
for each fixed \(u\), its derivative with
respect to \(u\) is also even in \(t\).
Hence
\[
D_u\mathcal S(u_0,t)[v]=Lv+O(t^2).
\]
Combining the preceding expansions yields
\[
\begin{aligned}
S_{tB_*}(u_t)
&=4-t^2Q+t^2Lv+O(t^4)\\
&=4+t^2(Lv-Q)+O(t^4).
\end{aligned}
\]

Since \(u_t\) satisfies the generalized
extremal equation
\[
S_{tB_*}(u_t)=4
\]
for all sufficiently small \(t\),
comparison of the coefficients of \(t^2\)
implies
\begin{equation}\label{eq:LvQ}
Lv=Q.
\end{equation}

Thus the leading correction \(v\) is governed
by a linear elliptic equation: \(L\) describes
the response of the ordinary scalar curvature
to variations of the symplectic potential,
whereas \(Q\) records the second-order effect
of the mixed parameter \(B_t\).
We compute these quantities explicitly below
and show that \(v\) has nonvanishing mixed
second derivatives.

\subsubsection{The product Fubini--Study structure}

Consider the standard labeled square
\[
P=[-1,1]\times[-1,1],
\]
with defining affine functions
\[
L_1=1+x,\quad L_2=1-x,\quad
L_3=1+y,\quad L_4=1-y.
\]
The corresponding Guillemin potential is
\begin{equation}\label{eq:base}
u_0(x,y)=f(x)+f(y),
\end{equation}
where
\[
f(s)=\frac12\left[
(1+s)\log(1+s)+(1-s)\log(1-s)
\right].
\]
This is the symplectic potential of the product
Fubini--Study metric, with the normalization
determined by the labeled polytope.

Set
\[
p(x)=1-x^2,\qquad q(y)=1-y^2.
\]
Since
\[
f''(s)=\frac1{1-s^2},
\]
we have
\[
H_0:=\Hess u_0=
\begin{pmatrix}
p^{-1}&0\\
0&q^{-1}
\end{pmatrix},
\qquad
H_0^{-1}=
\begin{pmatrix}
p&0\\
0&q
\end{pmatrix}.
\]
The ordinary Abreu formula gives
\[
S_0(u_0)=-p''(x)-q''(y)=4.
\]
Thus \(u_0\) is a constant-scalar-curvature
K\"ahler potential. The extremal affine function
of the labeled square is
\[
\ell_P=4.
\]

We introduce the family of skew-symmetric matrices
\[
B_t=tJ_*,
\qquad
J_*=
\begin{pmatrix}
0&1\\
-1&0
\end{pmatrix}.
\]
For the undeformed potential \(u_0\), we obtain
\[
H_0+\sqrt{-1} B_t=
\begin{pmatrix}
p^{-1}&t \sqrt{-1} \\
-t \sqrt{-1} &q^{-1}
\end{pmatrix},
\]
and hence
\[
\det(H_0+ B_t \sqrt{-1})
=\frac1{pq}-t^2.
\]
Since \(0<p,q\leq1\) in the interior of \(P\),
the Hermitian matrix \(H_0+\sqrt{-1} B_t\) is positive
definite for \(|t| \ll 1\). The corresponding
facewise admissibility conditions also hold.
Thus \(u_0\) remains admissible for sufficiently
small \(t\), although it does not remain
extremal when \(t\neq0\).

By Proposition \ref{prop:ift}, there exists
\(\varepsilon>0\) and a smooth family
\[
u_t\in\Sc_{B_t}(P),
\qquad |t|<\varepsilon,
\]
satisfying
\begin{equation}\label{eq:family-equation}
S_{B_t}(u_t)=4,
\qquad u_{t=0}=u_0.
\end{equation}
We choose the affine normalization
\[
u_t-u_0\perp\Aff(P)
\]
with respect to \(L^2(P,dx\,dy)\). Also we have

\begin{equation}\label{eq:even-expansion-42}
u_t=u_0+t^2v+O(t^4),
\end{equation}
where \(v\) is a smooth function on the closed
square. The remainder is understood in each
fixed \(C^k\)-norm for the smooth differences
of symplectic potentials.

\subsubsection{The second-order deformation equation}

We next determine the coefficient \(v\) in
\eqref{eq:even-expansion-42}. Note that
\begin{equation}\label{eq:exactinverse}
X_{B_t}(u_0)
=
\Ree(H_0+\sqrt{-1} B_t)^{-1}
=
\frac1{1-t^2pq}
\begin{pmatrix}
p&0\\
0&q
\end{pmatrix}.
\end{equation}
It follows that
\[
S_{B_t}(u_0)
=
-\partial_x^2\left(\frac{p}{1-t^2pq}\right)
-\partial_y^2\left(\frac{q}{1-t^2pq}\right).
\]
Expanding at \(t=0\), we obtain
\[
S_{B_t}(u_0)=4-t^2Q+O(t^4),
\]
where
\begin{equation}\label{eq:Q}
\begin{aligned}
Q(x,y)
&=(p^2)''(x)q(y)+p(x)(q^2)''(y)\\
&=-8+16x^2+16y^2-24x^2y^2.
\end{aligned}
\end{equation}

In particular, \(Q\) is not affine-linear.
Let \(L\) denote the linearization of
ordinary Abreu scalar curvature at \(u_0\).
Since
\[
\delta(H^{-1})
=-H_0^{-1}(\Hess v)H_0^{-1},
\]
we obtain
\begin{equation}\label{eq:L}
Lv=
\partial_x^2(p^2v_{xx})
+2\partial_x\partial_y(pqv_{xy})
+\partial_y^2(q^2v_{yy}).
\end{equation}
Substituting
\(u_t=u_0+t^2v+O(t^4)\)
into \eqref{eq:family-equation} and comparing
the coefficients of \(t^2\) gives
\begin{equation}
Lv=Q.
\end{equation}

\subsubsection{Explicit solution of the linearized equation}

Let
\[
P_2(s)=\frac12(3s^2-1)
\]
be the second Legendre polynomial.
The expression \eqref{eq:Q} can be rewritten as
\[
Q=
\frac{16}{3}\bigl(P_2(x)+P_2(y)\bigr)
-\frac{32}{3}P_2(x)P_2(y).
\]

We use the elementary identities
\[
\bigl((1-s^2)P_2'(s)\bigr)'=-6P_2(s),
\]
and
\[
\bigl((1-s^2)^2P_2''(s)\bigr)''
=24P_2(s).
\]
Applying \eqref{eq:L}, we find
\[
L(P_2(x))=24P_2(x),
\qquad
L(P_2(y))=24P_2(y).
\]
Moreover,
\[
\begin{aligned}
L(P_2(x)P_2(y))
&=(24+24+72)P_2(x)P_2(y)\\
&=120P_2(x)P_2(y).
\end{aligned}
\]

It follows that
\begin{equation}\label{eq:v}
v(x,y)=
\frac29\bigl(P_2(x)+P_2(y)\bigr)
-\frac4{45}P_2(x)P_2(y)
\end{equation}
satisfies \(Lv=Q\).

This solution is orthogonal to
\(\Aff(P)=\operatorname{span}\{1,x,y\}\).
It is also unique under this normalization.
Indeed, integration by parts gives
\begin{equation}\label{eq:Lpositive}
\int_P wLw\,dx\,dy
=
\int_P
\left(
p^2w_{xx}^2
+2pq\,w_{xy}^2
+q^2w_{yy}^2
\right)dx\,dy.
\end{equation}
The boundary terms vanish because \(p\) and
\(q\) vanish at the corresponding endpoints.
Consequently, \(Lw=0\) implies
\(\Hess w=0\) in the interior of \(P\), so
\(w\) is affine-linear.

We have therefore obtained the following
second-order expansion:
\begin{equation}\label{eq:expansion}
\begin{aligned}
u_t(x,y)
=u_0(x,y)+t^2\bigg[
&\frac29\bigl(P_2(x)+P_2(y)\bigr)\\
&-\frac4{45}P_2(x)P_2(y)
\bigg]+O(t^4).
\end{aligned}
\end{equation}

The crucial feature of this expansion is
the mixed term
\[
-\frac4{45}t^2P_2(x)P_2(y).
\]
Differentiating, we obtain
\begin{equation}\label{eq:mixedderivative}
\partial_x\partial_yu_t
=
-\frac45t^2xy+O(t^4).
\end{equation}
At the point \((x,y)=(1/2,1/2)\),
\[
\partial_x\partial_yu_t(1/2,1/2)
=
-\frac15t^2+O(t^4),
\]
which is nonzero for all sufficiently
small \(t\neq0\).

Since the mixed Hessian of any potential
of the form
\[
u(x,y)=u_1(x)+u_2(y)+a(x,y),
\qquad a\in\Aff(P),
\]
vanishes identically, this proves that
the deformed potentials \(u_t\) do not split.

We summarize the construction in the
following theorem.

\begin{theorem}[Perturbative non-splitting]
\label{thm:mixed}
Let \(P=[-1,1]^2\) be the standard labeled
moment polytope of
\(\mathbb{CP}^1\times\mathbb{CP}^1\), and let
\[
B_t=t
\begin{pmatrix}
0&1\\
-1&0
\end{pmatrix}.
\]
For sufficiently small \(|t|\), there exists
a smooth family of admissible toric generalized
constant-scalar-curvature K\"ahler potentials
\(u_t\in\Sc_{B_t}(P)\) satisfying
\[
u_{t=0}=u_0,\qquad S_{B_t}(u_t)=4.
\]
With the affine normalization described above,
this family has the expansion
\eqref{eq:expansion}. In particular,
\[
\Hess_{xy}u_t\not\equiv0
\]
for every sufficiently small nonzero \(t\).
Thus the associated symplectic potentials
do not split with respect to the prescribed
product decomposition.
\end{theorem}

\subsubsection{An exact obstruction to splitting}

The preceding calculation proves non-splitting
for sufficiently small parameters by exhibiting
a nonzero mixed term in the perturbation
expansion. We now prove a stronger statement:
for any nonzero mixed parameter, no splitting
potential can solve the generalized cscK
equation, provided it satisfies the admissibility
conditions.

\begin{proposition}[Exact non-splitting obstruction]
\label{prop:exactobstruction}
Let
\[
B_t=t
\begin{pmatrix}
0&1\\
-1&0
\end{pmatrix},
\qquad t\neq0,
\]
on the standard labeled square \(P=[-1,1]^2\).
There is no admissible symplectic potential
of the form
\[
u(x,y)=f_1(x)+f_2(y)
\]
modulo affine-linear functions satisfying
\[
S_{B_t}(u)=4.
\]
\end{proposition}

\begin{proof}
Suppose that such a potential exists. Set
\[
a(x)=\frac1{f_1''(x)},
\qquad
b(y)=\frac1{f_2''(y)}.
\]
The Guillemin boundary conditions imply
that \(a\) and \(b\) extend smoothly to
the closed intervals, are positive in
their interiors, and satisfy
\[
\begin{aligned}
a(\pm1)&=0,&
a'(-1)&=2,&a'(1)&=-2,\\
b(\pm1)&=0,&
b'(-1)&=2,&b'(1)&=-2.
\end{aligned}
\]
Indeed, for example,
\[
f_1''(x)=\frac1{1-x^2}+h_1''(x),
\]
where \(h_1\) is smooth on \([-1,1]\).
Consequently,
\[
a(x)=
\frac{1-x^2}
{1+(1-x^2)h_1''(x)},
\]
which gives the claimed boundary behavior.

By admissibility,
\[
1-t^2a(x)b(y)>0
\]
in the interior. The real part of the inverse
Hermitian matrix is
\[
X_{B_t}(u)=
\frac1{1-t^2ab}
\begin{pmatrix}
a&0\\
0&b
\end{pmatrix}.
\]
Hence the generalized scalar curvature
equation takes the form
\[
-\partial_x^2
\left(\frac{a}{1-t^2ab}\right)
-\partial_y^2
\left(\frac{b}{1-t^2ab}\right)
=4.
\]

We evaluate this identity on the edge
\(y=1\). Since \(b(1)=0\), the first term
reduces to \(-a''(x)\). Moreover,
\[
\left.
\partial_y^2
\left(\frac{b}{1-t^2ab}\right)
\right|_{y=1}
=
b''(1)+2t^2a(x)(b'(1))^2.
\]
Using \(b'(1)=-2\), we obtain
\begin{equation}\label{eq:edgea}
a''(x)+8t^2a(x)=-4-b''(1).
\end{equation}
Similarly, evaluating the equation
on the edge \(x=1\) gives
\begin{equation}\label{eq:edgeb}
b''(y)+8t^2b(y)=-4-a''(1).
\end{equation}

Integrating \eqref{eq:edgea} over
\([-1,1]\), and using
\[
a'(1)-a'(-1)=-4,
\]
we find
\[
b''(1)
=-2-4t^2\int_{-1}^1a(x)\,dx.
\]
Similarly,
\[
a''(1)
=-2-4t^2\int_{-1}^1b(y)\,dy.
\]

On the other hand, evaluating
\eqref{eq:edgea} at \(x=1\), where
\(a(1)=0\), yields
\[
a''(1)+b''(1)=-4.
\]
Combining these three identities gives
\[
t^2
\left(
\int_{-1}^1a(x)\,dx+
\int_{-1}^1b(y)\,dy
\right)=0.
\]
Since \(a,b>0\) in the interior and
\(t\neq0\), this is impossible.
The desired contradiction proves
the proposition.
\end{proof}

\begin{remark}
Theorem \ref{thm:mixed} produces exact
generalized cscK structures for all
sufficiently small nonzero \(t\), while
\eqref{eq:expansion} describes only their
leading deformation terms. The proposition
shows that non-splitting holds independently
of this asymptotic expansion, for every
nonzero parameter at which an admissible
cscGK solution exists.

The construction does not supply a closed
formula for \(u_t\) or a numerical lower
bound for the existence interval. Moreover,
varying \(B_t\) changes the holomorphic
Poisson tensor, so the resulting structures
generally belong to different generalized
K\"ahler classes.
\end{remark}

\printbibliography

@article{ASU,
  author = {Apostolov, Vestislav and Streets, Jeffrey and Ustinovskiy, Yury},
  title = {The {Riemannian} and symplectic geometry of the space of generalized {K\"ahler} structures},
  year = {2023},
  eprint = {2302.07314},
  archivePrefix = {arXiv},
  primaryClass = {math.DG}
}

@article{Gualtieri,
  author = {Gualtieri, Marco},
  title = {Generalized {K\"ahler} geometry},
  journal = {Communications in Mathematical Physics},
  volume = {331},
  pages = {297--331},
  year = {2014}
}

@article{LS,
  author = {LeBrun, Claude and Simanca, Santiago R.},
  title = {Extremal {K\"ahler} metrics and complex deformation theory},
  journal = {Geometric and Functional Analysis},
  volume = {4},
  pages = {298--336},
  year = {1994}
}

@article{ApostolovHuang,
  author  = {Apostolov, Vestislav and Huang, Hongnian},
  title   = {A Splitting Theorem for Extremal Kähler Metrics},
  journal = {Journal of Geometric Analysis},
  volume  = {25},
  number  = {1},
  pages   = {149--170},
  year    = {2015}
}

@article{HuangSplitting,
  author  = {Huang, Hongnian},
  title   = {A Splitting Theorem on Toric Manifolds},
  journal = {Mathematical Research Letters},
  volume  = {20},
  number  = {2},
  pages   = {273--278},
  year    = {2013}
}

@article{GatesHullRocek,
  author  = {Gates, S. James, Jr. and Hull, C. M.
             and Ro{\v{c}}ek, M.},
  title   = {Twisted Multiplets and New Supersymmetric
             Non-linear {$\sigma$}-Models},
  journal = {Nuclear Physics B},
  volume  = {248},
  number  = {1},
  pages   = {157--186},
  year    = {1984}
}

@article{Hitchin2003,
  author  = {Hitchin, Nigel},
  title   = {Generalized Calabi--Yau Manifolds},
  journal = {The Quarterly Journal of Mathematics},
  volume  = {54},
  number  = {3},
  pages   = {281--308},
  year    = {2003}
}

@phdthesis{GualtieriThesis,
  author  = {Gualtieri, Marco},
  title   = {Generalized Complex Geometry},
  school  = {University of Oxford},
  year    = {2004},
  eprint  = {math/0401221},
  archivePrefix = {arXiv}
}

@article{HitchinPoisson,
  author  = {Hitchin, Nigel},
  title   = {Instantons, Poisson Structures and
             Generalized K{\"a}hler Geometry},
  journal = {Communications in Mathematical Physics},
  volume  = {265},
  number  = {1},
  pages   = {131--164},
  year    = {2006}
}

@article{GibsonStreets,
  author  = {Gibson, Matthew and Streets, Jeffrey},
  title   = {Deformation Classes in Generalized
             K{\"a}hler Geometry},
  journal = {Complex Manifolds},
  volume  = {7},
  pages   = {241--256},
  year    = {2020},
  doi     = {10.1515/coma-2020-0101}
}

@incollection{GualtieriHamiltonian,
  author    = {Gualtieri, Marco},
  title     = {Generalized K{\"a}hler Metrics from
               Hamiltonian Deformations},
  booktitle = {Geometry and Physics: Volume II},
  year      = {2018}
}

\end{document}